\documentstyle[amssymb]{article}
          
\begin{document}
        \def\far{Funkt. Anal. i Prilozh.\ }
            \def\fae{Funct. Anal. Appl.\ }
                    \def\rrr{ (Russian)}

                               \def\eee{;  English translation:}

\def\square{\,\hbox{\vrule\vbox to 5pt
    {\hrule width 5pt \vfill\hrule}\vrule}\;}

\def\R{\mathbb R}
\def\C{\mathbb C}
\def\H{\mathbb H}
\def\K{\mathbb K}
\def\T{{\mathbb T}}

\title{Boundary values of holomorphic functions and
   some spectral problems  for
         unitary repesentations}
      \author{Yurii A. Neretin\thanks{supported by RFRF (grant 95-01-00814)
           and Russian program of support of scientific schools
           (grant 96-01-96249)}}
                                 \maketitle

Moscow State Institute of Electronics and Mathematics, Bolshoi Triohsvyatitelskii, per. 3/12, Moscow-109028, Russia

Max-Planck-Institut f\"ur Mathematik, Gottfried-Claren-Str. 26, 53225 Bonn, Germany

neretin \@ mpim-bonn.mpg.de, neretin \@ main.mccme.rssi.ru
\vspace{22pt}

Consider a (not necessarily irreducible) unitary representation of a semisimple Lie
group $G$ in a Hilbert space $H$ (for instance we can consider an action
of the group $G$ in some function space on some homogeneous space,
a tensor product of irreducible unitary representations, a restriction of
an irreducible unitary representation to a subgroup etc.). There arises
a natural problem: to decompose our representation into irreducible
representations, in other words, to find the spectrum of the representation.
Problems of this type were widely investigated for the last 50 years (since [Kre])
and by now many such spectral problems have been solved completely or
partially.

In several simple cases the spectra are continuous and more
or less `uniform'.  But sometimes such spectra contain
strange discrete increments which means that there exist
minimal $G$-invariant subspaces in the Hilbert space $H$.
First examples of discrete increments
were observed in [Nai], [Puk], [H-Ch], [Mol1]
(1961 - 1966) and now there exists a large
literature devoted to discrete spectra, see [Boy], [Ism],
[Mol3-4], [Str], [Far], [Pat], [F-J1], [F-J2], [Sch], [Kob1-2],
[RSW], [Tsu], [How], [Ada], [Li], [BO] etc.

Quite often `discrete increments'  are very singular (`exotic')
representations and it is very difficult to construct such unitary
representations in some other way. Even Harish-Chandra's discrete
series still cannot be considered handy objects. The remarkable
paper [F-J1] (1980, see also [F-J2]) contained a serious generalization
of Harish-Chandra's discrete series and this paper started a systematic
search of `exotic' unitary representations in discrete spectra.

In [Ner1] (1986), [Ols3],[Ols4], [NO], [Ner2]
 it was observed that very often
discrete increments to spectra are related to some functional analytic
phenomena, namely to so-called `trace theorems' (i.e theorems
about restrictions of discontinuous functions to submanifolds).
On various analytic trace theorems see [RS], IX.3,IX.9,
[ESh], 2.3, [Beu], [Bar],5.12, [NR], [Rud], [Vla], [VS], [Ner2].

This paper is some kind of addendum to papers [NO], [Ner2]
of Olshanskii and myself. Its purpose  is to discuss some open
problems related to applications of trace theorems in
noncommutative harmonic analysis.

 I am grateful to G.\ I.\ Olshanskii for our col\-lab\-ora\-tion. I also thank
V.\ F.\ Molchanov, H.\ Schlichtkrull, G.\ Zuckerman, B.\ {\O}rsted,
M.\ Flensted-Jensen, T.\ Kobayshi,\ A.\ G.\ Sergeev, G.\ \'Olafsson,
B.\ Kostant, A.\ Dvorsky, R.\ S.\ Ismagilov, R.\ Howe, V.\ M.\ Gichev,
V.\ V.\ Lebedev for discussions, comments and references. These notes are
mainly based on my lectures [Ner5] delivered at
the school `Analysis on homogeneous spaces'
(Tambov, august 1996) and I also thank the organizers
of that school.

\section{ Tensor product of two representations of \\
 $SL(2,\R)$   of the complementary series}

{\bf 1.1. Complementary series.}
   Consider  the circle $\T^1$ with the
 coordinates $z=e^{i\phi}$ (where $| z| =1
   , \phi\in [0,2\pi]$).  Let $-1<s<1,s\neq 0$.
   Further consider the scalar product
   \begin{equation}\label{dop}
         <f_1,f_2>=\int_0^{2\pi}\int_0^{2\pi}
      {
          \frac{f_1(\phi_1)\overline{f_2(\phi_2)}d\phi_1d\phi_2}
            {|\sin (\phi_1-\phi_2)/2|^{(1+s)/2}}
        }
           \end{equation}
      in the space of $C^\infty$-functions on the circle.
      The functions $e^{in\phi}$ are pairwise orthogonal with respect
      to this scalar product and
         \begin{equation}\label{gamma}
      <e^{in\phi},e^{in\phi}>=
      \frac{\Gamma(n+(1+s)/2)}{\Gamma(n+(1-s)/2)}
      \end{equation}
If  $s>0$, then the integral (\ref{dop}) diverges. Nevertheless,
the formula (\ref{gamma})  makes sense and we can consider (\ref{gamma})
as the definition of the divergent integral  (\ref{dop}).
 Denote by ${\cal H}_s$ the completion  of the
   space $C^\infty (\T^1)$    with respect to the scalar product (\ref{dop}).
    We have
   $$
    f(\phi)=\sum c_n e^{in\phi}\in {\cal H}_s  \Leftrightarrow
    \sum  n^{s}|c_n|^2 <\infty
    $$
    i.e. ${\cal H}_s$ is a {\it Sobolev space} on the circle $\T^1$.

     We realize the group $SL(2,\R)$ as the group of matrices
    $g=\left( \begin{array}{cc}
           a&b\\ \overline b &\overline a\end{array}\right)   $
    satisfying the condition $|a|^2-|b|^2=1$.
Define a representation of the group $SL(2,\R)$ in the space
${\cal H}_s$ by the formula
   $$
           T_s\left( \begin{array}{cc}
            a&b\\ \overline b &\overline a\end{array}\right) f(z)
           =f\left(
           \frac {az+b}{\overline bz+\overline a}
                 \right)|\overline b z+ \overline a|^{s-1}.
                 $$
           The representations $T_s$ are unitary
in the Hilbert space ${\cal H}_s$ (this fact is simple but
           it is not trivial). They are called the
{\it complementary series representations of $SL(2,\R)$.}
Note that the representations $T_s$ and $T_{-s}$ are equivalent.

  Consider the tensor product  $T_{s_1}\otimes T_{s_2}$
of two complementary series representations.
The space  ${\cal H}_{s_1}\otimes {\cal H}_{s_2}$ is a
space of functions on the two-dimensional torus $\T^2:
 z_1=e^{i\phi_1},z_2=e^{i\phi_2}$
 equipped with the scalar product
\begin{eqnarray*}
& &   <f_1,f_2> \\
&=& \int_0^{2\pi}\int_0^{2\pi}\int_0^{2\pi}\int_0^{2\pi}
  \frac {f_1(\phi_1,\phi_2) \overline{
  f_2(\psi_1,\psi_2)}d\phi_1d\phi_2d\psi_1
  d\psi_2}{
|\sin (\phi_1-\psi_1)/2|^{(1+s_1)/2}|\sin (\phi_2-\psi_2)/2|^{(1+s_2)/2}}
\end{eqnarray*}
 The group $SL_2(\R)$ acts in this space via the formula
   $$   (T_{s_1}\otimes T_{s_2})
   \left( \begin{array}{cc}
    a&b\\ \overline b &\overline a\end{array}\right) f(z_1,z_2)=
   f\left(
   \frac {az_1+b}{\overline bz_1+\overline a},
   \frac {az_2+b}{\overline bz_2+\overline a}\right)
   |\overline bz_1 + \overline a|^{s_1-1}|\overline bz_2+ \overline a|^{s_2-1}
   $$

{\bf 1.2. Restriction to the diagonal.}

Denote by $\Delta$ the diagonal $\phi_1=\phi_2$ of the torus $\T^2$.

{\bf Theorem 1.1.} {\it Let $s_1+s_2>1$ . Then the restriction
operator restricting a function $f\in C^\infty(\T^2)$ to the diagonal
$\Delta$ extends to a bounded operator $R$
from the space  ${\cal H}_{s_1}\otimes {\cal H}_{s_2}$
to the space $L^2(\Delta)$. }

  We emphasize that the functions $f\in
    {\cal H}_{s_1}\otimes {\cal H}_{s_2}$
 may be discontinuous and hence a function $f\in {\cal H}_s$
 has no value in a individual point of the torus.
  Nevertheless, the restriction operator $R$ is well defined.

{ Proof}.   The space  ${\cal H}_{s_1}\otimes {\cal H}_{s_2}$
 consists of functions
 $$f(\phi_1,\phi_2)=\sum c_{n,m} e^{in\phi_1}e^{im\phi_2}$$
 satisfying the condition
 $$
 \sum (1+|n|)^{s_1}(1+|m|)^{s_2} |c_{n,m}|^2<\infty
 $$
 The restriction of $f$ to the diagonal $\Delta$
 has the form
 $$
 \sum_k (\sum_{m+n=k} c_{m,n} )e^{ik\phi}
 $$
 and the convergence of the series $\sum c_{m,n}$ easily
 can be deduced from the Cauchy-Bunyakovsky inequality.
$\square$

 The restriction operator is a bounded $SL(2,\R)$-interwining
 operator
 from   $T_{s_1}\otimes T_{s_2}$ to $T_{s_1+s_2-1}$. The representation
  $T_{s_1}\otimes T_{s_2}$ is unitary and hence $T_{s_1+s_2-1}$
 is a subrepresentation of  $T_{s_1}\otimes T_{s_2}$.

 {\bf 1.3. Dual picture: distributions supported on the diagonal.}
 Now let $s_1+s_2<-1$ (recall that $T_{-s}\simeq T_s$).
 It turns out  that in this case  the Hilbert
 space  ${\cal H}_{s_1}\otimes {\cal H}_{s_2}$
 contains  distributions supported on the diagonal $\Delta$.
 Obviously the space of distributions supported
 on the diagonal is an invariant subspace.
 It is easy to show that the representation of $SL(2,\R)$
 in this subspace is equivalent to $T_{s_1+s_2+1}$
 and we again obtain an embedding of a  complementary series
 repesentation into the tensor product.

 {\bf 1.4. Separation of a discrete increment.}
 Let $s_1=s_2=s<-1/2$. Consider the operator
 $$
 Jf(\phi_1,\phi_2)= \left|\sin\left(\frac{\phi_1-\phi_2}
{2}\right) \right|^{s}f(\phi_1,\phi_2)
 $$
 Then the operator $J$ interwines the representation
  $T_{s}\otimes T_{s}$ and the representation
 $$
 (T_{0}\otimes T_0) \left( \begin{array}{cc}
       a&b\\ \overline b &\overline a\end{array}\right) f(z_1,z_2)=
   f(
      \frac {az_1+b}{\overline bz_1+\overline a},
      \frac {az_2+b}{\overline bz_2+\overline a})
      |\overline bz_1 + \overline a|^{-1}
 |\overline bz_2+ \overline a|^{-1}.
 $$
 The last representation is the
  standard representation of $SL(2,\R)$
 in $L^2(\T^2)$. The kernel of $J$ consists of distributions
 supported on the diagonal  and the image is dense in $L^2$.
  Hence  the representation $T_s\otimes T_s$
 is equivalent to the direct sum of the representation $T_{1+2s}$
 and the representation of $SL(2,\R)$ in $L^2$ on the torus.

{\bf Remark.}  The operator $J$ is
a operator of multiplication on nonsmooth function in space of
distributions. This operator is
 unbounded and hence it is nessesary to be careful, see [NO], \S 7
 for details and for the case $s_1\neq s_2$.

{\bf Remark.}  The existence of the embedding of $T_{s_1+s_2-1}$
 into  $T_{s_1}\otimes T_{s_2}$ was obtained   by Pukanzky [Puk]
 in 1961, and an analogous result for $SL(2,{\C})$ was obtained by
Naimark  [Nai] (1961). The explanation of this phenomenon using
the restriction  operator wasn't observed before  [Ner1] (1986).

 \section{\bf Hilbert spaces defined by  positive defined kernels.}

 {\bf 2.1. Positive definite kernels.}
 Let $\Xi$ be a set. A function $K(x,y)$ on $\Xi\times\Xi$
 is called a {\it positive definite kernel} if for each
 finite collection of points $x_1,\dots,x_n\in \Xi$ we have
 $$
 \det\left(
 \begin{array}{cccc}
 K(x_1,x_1) & K(x_1,x_2)&\dots  & K(x_1,x_n)\\
 K(x_2,x_1)& K(x_2,x_2)&\dots & K(x_2,x_n)\\
 \vdots &\vdots &\ddots& \vdots\\
 K(x_n,x_1)&K(x_n,x_2) &\dots & K(x_n,x_n)
 \end{array}\right)
 \ge 0
 $$

 To each positive definite kernel $K(x,y)$ we associate a Hilbert space
 $H$ and fixed system of vectors $\Psi_x\in H$ enumerated by $x\in\Xi$
 such that

 1. The linear span of the system $\Psi_x$ is dense in $H$

 2. For all $x,y\in \Xi$ we have
 \begin{equation}\label{posit}
 <\Psi_y,\Psi_x>_H=K(x,y)
 \end{equation}

 The system $\Psi_x$ is called a {\it supercomplete basis}
or an {\it overfull system}.

Let $h\in H$. We consider the function   $f_h(x)=Jh(x)$ on $\Xi$
given by formula
 \begin{equation}\label{trans}
 JF(x)=f_h(x)=<h,\Psi_x>
 \end{equation}
  where $x\in\Xi$.

 {\it We will identify the space $H$ and
the image of the operator $J$.} In other words {\it we
consider the Hilbert space $H$ as a space of functions
on $\Xi$.}

 {\bf 2.2. Hilbert spaces of holomorphic functions.}
 Let $\Omega$ be an open domain in ${\C}^n$.
 A function $K(z,u)$ on $\Omega\times\Omega$ is called
 a {\it reproducing kernel} on $\Omega$ if

 1. $K(z,u)$ is a positive definite kernel

 2. The function $K(z,u)$ is holomorphic with respect to the
variable $z$  and antiholomorphic with respect to the
variable $u$.

 Consider  the associated  Hilbert space $H$ and let us realize
 it as a space of functions on $\Omega$ by the formula (\ref{trans}).
 It is easy to see that the functions $f_h$ are holomorphic, i.e. the
 Hilbert space $H$ can be considered as a Hilbert space of holomorphic
 functions on $\Omega$.

 {\bf 2.3. Example:  Berezin spaces on Cartan domains.}
Let $p\le q$.  Denote by $B_{p,q}$ the space of all $p\times q$-matrices $z$ satisfying
 the condition $\|z\|<1$, i.e $B_{p,q}$ is a {\it Cartan domain
 of the first type}.

  The pseudounitary  group  $U(p,q)$ consists of $(p+q)\times (p+q)$-matrices
   $
   g=\left(\begin{array}{cc} a&b\\ c&d \end{array}\right)
   $
   satisfying the condition
   \begin{equation}\label{upq}
   g\cdot\left(\begin{array}{cc} 1&0\\0&-1\end{array}\right)\cdot g^\ast
   =\left(\begin{array}{cc} 1&0\\0&-1\end{array}\right)
 \end{equation}
   The group $U(p,q)$ acts on $B_{p,q}$ by the
{\it fractional-linear} transformations
  \begin{equation}\label{fra}
   z\mapsto z^{[g]}:=(a+zc)^{-1}(b+zd)
  \end{equation}
   The stabilizer of the point $z=0$ consists of matrices having the form
   $\left(\begin{array}{cc} a&0\\ 0&
                                d\end{array}\right),$
   where $a\in U(p),d\in U(q)$. Hence $B_{p,q}$ is the
Riemannian symmetric space
   $$U(p,q)/(U(p)\times U(q))$$

{\bf Theorem 2.1.} {\rm(see [Ber1], see also [Wal])}
{\it Let  $s=0,1,2,\dots,p-1$ or $s>p-1$.    Then the kernel
   $$
   K_s(z,u)={\det}^{-s}(1-zu^\ast)
   $$
   is positive definite on $B_{p,q}$.}

   Denote by $V_s$ the Hilbert space of holomorphic functions on $B_{p,q}$
   defined by the kernel $K_s(z,u)$ (for a discussion of this space
 see for instance [Ber1], [NO]).
    The group $U(p,q)$ acts on $V_s$ by the unitary operators
  \begin{equation}\label{act}
  T_s\left(\begin{array}{cc} a&b\\ c&d\end{array}\right)f(z)
   =f(z^{[g]}){\det}^{-s}(a+zc),
 \end{equation}
   where $z^{[g]}$ is given by the formula (\ref{fra}).

 The representations $T_s$ are called  {\it spherical
  highest weight representations} of $U(p,q)$.

{\bf Remark.}  If $s$ is an integer, then $T_s$ is a representation
   of the group $U(p,q)$ itself. If $s$ is not an integer,
 then the function $\det^{-s}(\cdot)$ is multivalued and hence  $T_s$ is a
   representation of the universal covering group of $U(p,q)$
 (in other words, the representation $T_s$ is a
 projective representation of the group $U(p,q)$ ).

{\bf Remark.}  If $s>p-1$, then the space
 $V_s$ contains all polynomial functions
 on the Cartan domain $B_{p,q}$. If $s=0,1,\dots,p-1$, then all functions
 $f\in V_s$ satisfy the following system of partial  differential equations:
 $$
\left[ \det\left(\begin{array}{ccc}
 \frac{\partial}{ {\partial z_{i_1 j_1}}}&\dots&
 \frac{ \partial}{{\partial z_{i_\alpha j_1}}}
 \\ \vdots&\ddots&\vdots \\
 \frac{\partial}{{\partial z_{i_1 j_\alpha}}}&\dots&
 \frac{\partial}{{\partial z_{i_\alpha j_\alpha}}}
 \end{array}\right)
 \right] f(z)=0
 $$
 for all $\alpha=2,3,\dots,s+1$ and all $i_\mu$, $j_\nu$ satisfying
the inequalities $0\le i_\mu\le q$, $0\le j_\nu\le p$.

 {\bf 2.4. Cartan domains of second and third types.}
 {\it The Cartan domain $C_n$ of  third type} is the set of all $n\times n$
 symmetric matrices $z$ with complex coefficients satisfying the
 condition  $\|z\|<1$.
 The group $G$ of biholomorphic transformations (automorphisms)
 of the domain $C_n$  consists of transformations having the
form (\ref{fra}),  where a matrix
$g= \left(\begin{array}{cc}a&b\\ c&d \end{array}\right)$
 satisfies the condition
 \begin{equation}\label{sp}
 \left(\begin{array}{cc}a&b\\ c&d \end{array}\right)
 \left(\begin{array}{cc}1&0\\ 0&-1 \end{array}\right)
 \left(\begin{array}{cc}a&b\\ c&d \end{array}\right)^t=
 \left(\begin{array}{cc}1&0\\ 0&-1 \end{array}\right)
 ;\quad d=\overline a,\quad c= \overline b.
 \end{equation}
 It is easy to see that the group $G$ is isomorphic to the real symplectic
 group $Sp(2n,\mathbb R)$ and the domain $C_n$ is the Riemannian
symmetric space
 $$
 C_n=Sp(2n,{\mathbb R})/U(n).
 $$
 The kernel
 $$
 K_s(z,u):={\det}^{-s} (1-z\overline u)
 $$
 is positive definite if and only if
 $$
 \{s=0,1/2,\dots,(n-1)/2\}\ {\mbox or}\  \{ s>(n-1)/2\}.
 $$

   The {\it Cartan domain $D_n$} of the second type is the set of all
   $n\times n$ skew-symmetric matrices $z$ satisfying the condition
   $\|z\|<1$. The group $G$ of automorphisms of the domain $D_n$
   is the group of transformations (\ref{fra}),  where the matrix
 $g=  \left(\begin{array}{cc}a&b\\ c&d \end{array}\right)$
 satisfies the condition
 \begin{equation}\label{so}
 \left(\begin{array}{cc}a&b\\ c&d \end{array}\right)
 \left(\begin{array}{cc}0&1\\ 1&0 \end{array}\right)
 \left(\begin{array}{cc}a&b\\ c&d \end{array}\right)^t=
 \left(\begin{array}{cc}0&1\\ 1&0 \end{array}\right)
 ;\quad d=\overline a, \quad c=\overline b.
\end{equation}

{\bf Remark.}  Condition (\ref{so}) implies (\ref{upq}).

The group $G$ is isomorphic to the real classical
group $SO^\ast(2n)$ and the domain $D_n$ is the symmetric space
$$
D_n=SO^\ast(2n)/U(n).
$$
The kernel $K(z,u)=\det^{-s}(1-zu^\ast)$
on $D_n$ is positive definite if $s=0,1,2,\dots,\\ n-1$ or $s>n-1$.

{\it Spherical highest weight representations}
 of $Sp(2n,\mathbb R)$ and $SO^\ast(2n)$
act in Hilbert spaces associated to the kernels
$K_s$ by the formula (\ref{act}).

{\bf 2.5. Matrix-valued positive definite kernels.}
 Consider a set $\Xi$  and a finite dimensional complex euclidian space
   $\cal Y$. Denote by $GL(\cal Y)$ the group of invertible linear
   operators in $\cal Y$.     We say that a function
   $$
   L:\Xi\times\Xi \rightarrow GL({\C},\cal Y)
     $$
    is a {\it matrix-valued positive definite kernel}
(see for instance [HN]) if the function
     $$
     \widetilde L((z,\xi);(u,\eta))
    := <L(z,u)\xi , \eta>\ ;\quad (z,\xi),(u,\eta)\in  \Xi\times \cal Y
      $$
      is a positive definite kernel on $\Xi\times \cal Y$.

    Suppose we have a matrix-valued positive defined kernel on $\Xi$. Then
   there exists a unique Hilbert space $H$  and a map $\Psi:\Xi\times{\cal Y}
    \rightarrow H$ such that

   a) The map $\Psi$ is linear on each fiber $z\times {\cal Y}\subset
 \Xi \times \cal Y$.

      b) $<\Psi(z,\xi),\Psi(u,\eta)>_H=<L(z,u)\xi , \eta>_{\cal Y}$

   c) The image of the map $\Psi$ is total in $H$.

 Again there exists a canonical map $J$ from the Hilbert space $H$
  to the space of $\cal Y$-valued functions on $\Xi$. Let $h\in H$.
 Then the function $f_h=Jf$ is defined via the equation
 $$<f_h(z),\xi>_{\cal Y}=<h,\Psi(z,\xi)>_H.$$

{\bf 2.6. Highest weight representations of the groups}
$G=U(p,q)$, $Sp(2n,{\mathbb R})$, $SO^\ast (2n)$.

 Consider the case $G=U(p,q)$. Let $\rho=\rho_1\otimes \rho_2$
be an irreducible  finite-dimensional  holomorphic representation
of the group  $GL(p,{\C})\times GL(q,{\C})$
or of its universal covering. Further let
$W=W_1\otimes W_2$
be the space of the representation $\rho$.
Assume that the function
$$
K_\rho (z,u):=\rho_1(1-zu^\ast)\otimes \rho_2(1-u^\ast z)
$$
is a matrix-valued positive-definite kernel (see
[Ols2] for conditions of positive definiteness, see also
[NO]).
Consider the associated Hilbert space $H_\rho$ of
 $W_1\otimes W_2$-valued holomorphic functions
 on $B_{p,q}$. Then the {\it unitary highest
 weight representation} $T_\rho$ of the group
 $U(p,q)$ (or of its universal covering)
 acts on the space $H_\rho$ by the formula
 \begin{equation}  \label{hw}
 T_\rho
 \left( \begin{array}{cc}
 a&b\\ c&d \end{array}
 \right) F(z)=
 (\rho_1(a+zc)\otimes \rho_2(b-dz^{[g]}))
F(z^{[g]})
\end{equation}

Let $G=Sp(2n,\R)$ or $SO^\ast(2n)$.
Consider a finite dimensional holomorphic representation
$\rho$ of the group $GL(n,{\C})$ or of its universal
covering group. Assume  the matrix-valued kernel
$K_\rho(z,u)={\det}^{-s}(1-zu^\ast)$
to be positive-definite on the Cartan domain
$C_n$ or $D_n$, respectively. Denote by $H_\rho$ the
associated Hilbert space of holomorphic vector-valued functions.
   Then the {\it unitary highest weight representation}
   $T_\rho$ of the group $G$ is given by the formula
   $$
T_\rho
 \left( \begin{array}{cc}
 a&b\\ c&d \end{array}\right)
 f(z)=\rho (a+zc) F(z^{[g]}).
$$

{\bf 2.7. Bibliographical comments.}

a) There are also {\it Cartan domains of 4-th type}
$SO(n,2)/(SO(n)\times O(2))$ and two exceptional Cartan domains.
There exist associated highest weight representations.
The Cartan domains $SO(n,2)/(O(n)\times O(2))$ are a popular
object of complex analysis ( {\it tube of future }), see [VS].
These domains often are called  {\it Lie spheres or Lie balls},
but this term is a result of a philological mistake%
.
(Lie sphere is a oriented sphere, oriented plane or point in ${\R}^n$.
The space of all Lie spheres is a homogeneous space
$O(n+2)/(O(n)\times O(2))$ (see for instance [Cec]). The term
'Lie sphere' for a homogeneus space $O(n,2)/(SO(n)\times O(2))$
appeared in russian translation of famous chinese book [Hua].
It is quite strange to name this space by  'sphere' and it was transformed
to 'ball'.)

b) Conditions for the positive-definiteness of the standard
scalar-valued kernels on Cartan domains of types one through
four were discovered by Berezin [Ber] (with some additions
given by Gindikin [Gin]); Berezin and Gindikin
didn't consider exceptional domains; later
conditions for the positive-definiteness in the general
case (including exceptional)
were obtained by Wallach [Wal], Rossi - Vergne [RV] and Jackobsen [Jac],
proofs of these authors are essentially different).

c) Conditions for the positive-definiteness of the standard
kernels for $B_{p,q}$ were obtained in [Ols2].
For other domains this problem was discussed by many authors,
the final result was obtained in [EHW].

 \section{ Boundary values of holomorphic functions.}

 {\bf 3.1. $L^1$-limits.}

 Let  $\Omega \subset {\mathbb C}^n$ be an open domain, $\partial
 \Omega$ its boundary, and $\overline\Omega$ be the closure
  of $\Omega$. We say    that
  $\Omega $ is a {\it regular circle domain} if

  a) for all $z\in\Omega$ and $\lambda\in\mathbb C$ such that
$|\lambda| \le 1$   we have $\lambda z\in\Omega$

  b) for all $z\in\partial\Omega$ and $\lambda\in\mathbb C$ such that
   $|\lambda| <1$    we have $\lambda z\in\Omega$.

    Let $K(z,u)$ be a reproducing kernel
     in $\Omega$ satisfying the condition
   \begin{equation}\label{inv}
     K(e^{i\phi}z,e^{i\phi}u)=K(z,u)
   \end{equation}
    and $H$ be the Hilbert space of holomorphic functions
     associated to this kernel.

{\bf Remark.}  Cartan domains are regular circle domains
   and the Berezin kernels satisfy the condition (\ref{inv}).

{\bf Theorem 3.1.} {\rm(see [NO])} {\it Let $M\subset\partial\Omega$
   be a compact subset.
     Let $\mu$ be a measure supported on $M$. We make the
following assumptions

a) $K^\ast(z,u):=\lim_{c\rightarrow 1-0} K(cz,cu)$ exists almost
everywhere  on $M\times M$ with respect to the measure $\mu\times\mu$.

b) $K^\ast \in L^1(M\times M,\mu\times\mu)$
and  $\lim_{c\rightarrow 1-0}K(cz,cu)$
is dominated, i.e. there exists a function
$S(z,u)\in L^1(M\times M,\mu\times\mu$)
such that $|K(cz,cu)|<S(z,u)$ almost everywhere on $M\times M$ .

Then the restriction operator  restricting a function $f\in H$ to the set $M$
 is a well-defined operator
 $$
  H\rightarrow L^1(M,\mu).
  $$}

 This theorem doesn't cover all cases in which the restriction
operator exists (see Subsections 3.2-3.3 below). Hence
 the following natural problem arises.

 {\bf Problem 1}.{ Let $\Omega\subset
{\mathbb C}^n$ be an open domain. Let $K(z,u)$
 be a reproducing kernel and let $H$ be
the associated Hilbert space. Further let $M$ be
 a submanifold of the Shilov boundary of $\Omega$.
 Find conditions for the   existence of a
restriction operator from $H$ to some  Hilbert space
of functions on M.}

{\bf Remark.}  The case which is interesting for harmonic analysis is
the following. Let the domain $\Omega$ be a Cartan domain
$G/K$ and $H=V_s$ be a Berezin space.
 Consider a subgroup $Q\subset G$ and some $Q$-orbit
$M$ in the Shilov boundary of the domain $\Omega$. Suppose
the restriction operator restricting functions $f\in H_s$ to the manifold
$M$ exists. Then we obtain an interwining operator from the
Berezin space $H_s$ to some Hilbert space $H(M)$ of functions on $M$.
The representation of $Q$ in $H(M)$ is equivalent to a subrepresentation
of the representation of $Q$ in $H$.
 If the group  $Q$ is semisimple and the orbit $M$ is compact,
then the representation of $Q$ in $H(M)$ has purely discrete spectrum.

{\bf 3.2. Restrictions of functions of polynomial growth on a ball.}
 Denote by $B_q$ the unit ball
  $$
 |z_1|^2+\dots +|z_q|^2<1
  $$
in ${\mathbb C}^q$ and let $H^\infty(B_q)$ be the
Hardy space of bounded holomorphic functions on $B_q$.
Denote by $\langle\cdot,\cdot\rangle$ the standard
scalar product in ${\C}^n$.

 Consider a $C^1$-curve $\gamma(t)$
 in the sphere $\partial B_q=S^{2q-1}$.

{\bf Theorem 3.2} {\rm(see [NR],[Rud])}
{\it Let $\gamma(t)$ satisfy the condition

  \begin{equation}\label{nonta}
  \forall t: Im\langle\gamma(t),\gamma'(t)\rangle\neq 0
  \end{equation}

  Then for each $f\in H^\infty(B^q)$ the nontangential limit
 $f(z)$ as $z\rightarrow \gamma(t)$
exists almost everywhere on the curve $\gamma(t)$.}

{\bf Remark.}  Obviously the vector $\gamma'(t)$ is tangent to the sphere
   $S^{2q-1}$. Hence we have $Re\langle\gamma(t),\gamma'(t)\rangle=0$
 for all $t$. Fix $z\in S^{2q-1}$ and let $v\in {\C}^q$ be a vector.
The equation $\langle z,v\rangle=0$
 defines a subspace of codimension 1 in the  tangent space
 to the sphere $S^{2q-1}$ in the point $z\in S^{2q-1}\subset{\C}^n$.
 We obtain a contact distribution
 on the sphere (see [Rud], 5.4, 10.5). The condition (\ref{nonta})
  means that $\gamma(t) $ is transversal to this distribution.

 Denote by ${\cal D}'$ the space of all
 holomorphic {\it functions of polynomial
   grouth} in  $B^q$:
  $$f\in {\cal D}'\Leftrightarrow \exists N: \sup|f(z)|(1-|z|^2)^N<\infty $$

{\bf Remark.}  The space ${\cal D}'$ is a union
 ${\cal D}'=\cup_{s>0} V_s(B_q)$ of Berezin spaces on the ball, see
Subsection 2.3. The space $V_q(B_q)$ coincides with the
Hardy space $H^2(B_q)$.
The space $V_{q+1}(B_q)$ is the Bergman space (i.e the
intersection of $L^2(B_q)$ with the space of holomorphic functions).

{\bf Remark.}  It is well known (and more or less obvious)
 that each function
  $f\in {\cal D}'$ has a limit on the boundary in the sense of
  distributions (see [RS],IX.3 for a discussion
    of  theorems of this type and further references).

Denote by $R$ the restriction operator resctricting a
holomorphic function $f$  on the ball $B_q$ to the curve
$\gamma(t)$. Obviously the function $Rf$ is well-defined
if the function $f$ is continuous up to the boundary.

{\bf Theorem 3.3.} {\rm(see [Ner2])}
{\it Let $\gamma(t) $ be a $C^\infty$-smooth curve
in $\partial  B^q$ satisfying the condition (\ref{nonta}).
Then the restriction operator $R$
restricting a holomorphic function  to $\gamma(t)$
extends to a  bounded operator from ${\cal D}'$ to the
space of distributions on $\gamma(t)$.}

  {\bf 3.3. Restriction of functions of polynomial
growth on a polydisc.}

  Denote by $P^n$ the polydisc
  $|z_1|<1,\dots,|z_n|<1$. Let $\T^n$ be the
   torus $z_1=e^{i\phi_1},\dots,z_n=e^{i\phi_n}$.
Let $\gamma(t)=(\phi_1(t),\dots,\phi_n(t)) $
   be a {\it time-like} $C^\infty$-curve in $\T^n$ i.e.
  \begin{equation}
  \forall t: \phi_1'(t)>0,\dots,\phi'_n(t)>0
  \end{equation}
 Denote by ${\cal D}'$ the space of holomorphic
   functions of polynomial growth in $P^n$,
   i.e
   $$
   f\in {\cal D'}
\Leftrightarrow \exists N: \sup f(z) \prod_{k=1}^n (1-|z_k|)^N<\infty
   $$
    This condition is equivalent to the existence of an $L > 0$ with
   $$
   f(z_1,\dots,z_n)=\sum
c_{k_1,\dots,k_n}z_1^{k_1}\dots z_n^{k_n}\in{\cal D'}
\quad \hbox{ and } \quad
$$
$$ \sum \frac{|c_{k_1,\dots,k_n}|}{ [(k_1+1)\dots (k_n+1)]^L}<\infty
   $$

{\bf Theorem 3.4.} {\rm(see [Ner2])}
{\it  The restriction operator $R$  restricting a holomorphic function
  $f$ on $P^n$ to the curve $\gamma(t)$  extends to a  bounded operator
  from the space ${\cal D}'$ to the space
 of distributions on the curve $\gamma(t)$.}

 {\bf Example.} Consider a bidisc $P^2$ and
 the Hilbert space $H_{s_1,s_2} $ associated to the reproducing
 kernel
 $$
 (1-z_1\overline u_1)^{-s_1} (1-z_2\overline u_2)^{-s_2},
 $$
 where $s_1,s_2>0$. It is easy to show that
 ${\cal D}'=\cup_{s_1,s_2>0} H_{s_1,s_2}$.
  If $\gamma(t)$ has the form
 $(\gamma_1(t),\gamma_2(t))=(e^{it},e^{i\sigma})$
  where $\sigma=const$,
   then a restriction operator to $\gamma(t)$
   doesn't exist for all $s_1,s_2$. If $\gamma(t)$ satisfies
   the condition $\gamma'_1(t)>0,\gamma'_2(t)<0$ then the restriction
operator exists only in the case $s_1+s_2<1$.

{\bf 3.4. A conjecture.}
It seems to me that Theorems 3.3-3.4 are partial
cases of some general fact.
 Let $\Omega\subset{\mathbb C}^N$
be an open domain. Let $M$ be a submanifold in the Shilov boundary of
$\Omega$. Denote by $T_m$ the tangent
 space to $M$ in the point $m\in M$.
We identify $T_m$ with a
 linear submanifold in ${\mathbb C}^n$. Denote by $S_m \:=i\cdot T_m$ the
  linear submanifold which consists of the vectors
  $$
  i\cdot (v-m)+m,
  $$
  where $v\in T_m$ (and $i^2=-1$).
We say the submanifold $M$ is {\it time-like}
if  for each point $m\in M$ there
exists an open cone  $C_m\subset S_m$
with the vertex $m$ and an \\ $\epsilon$-neighborhood
 ${\cal U}_\epsilon(m)$ of $m$ such that
  $$
  C_m\cap \Omega \supset {\cal U}_\epsilon(m)\cap C_m.
  $$

  {\bf Conjecture.} Each holomorphic function
 of polynomial growth in
$\Omega$  has a restriction to a time-like manifold $M$
in the Shilov boundary in the sense of distributions.

    Theorems 3.3-3.4 are partial cases of this conjecture.
    It is also similiar to the standard
   facts on limits of functions of
  polynomial growth on the
 {\it whole} Shilov boundary mentioned above
   (see [RS], IX.3, [Vla], [VS]).
   Neverless I couldn't find this fact in the literature.

\section{  Positive definite kernels on
Riemannian noncompact symmetric spaces}

{\bf 4.1. Spherical kernel-representations of the groups $U(p,q)$.}
Consider the matrix ball $B_{p,q}$, see Subsection 2.3.
Let  $s=0,1,2,\dots,p-1$ or $s>p-1$.
 Consider the function
 $$
   L_s(z,u)=|\det (1-zu^\ast)|^{-2s},
   $$
   where $z,u\in B_{p,q}$.
With the notations of Subsection 2.3 we have
$$
L_s(z,u)=K_s(z,u)\overline{K_s(z,u)}
$$
and hence $L_s(z,u)$ is a positive definite
 kernel on $B_{p,q}$.
Denote by $H_s$ the Hilbert space
 associated with the kernel
$L_s(z,u)$. Using the formula (\ref{trans})
  we identify this space
with some space of real analytic functions on $B_{p,q}$.

   The group $U(p,q)$ acts on $H_s$ by
  the unitary operators
 $$
  A_s(g)f(z)=f(z^{[g]})|\det(a+zc)|^{-2s},
  $$
 where $z^{[g]}$ is defined by the formula (\ref{fra}).

  {\bf Problem 2.} {\it Decompose the representation $A_s$}

   We will call the representations $T_s$  {\it spherical
  kernel-representations} of the group $U(p,q)$.

Some partial cases of this question
were discussed in the end of the 70-ies,
  see [Ber2],[Rep],[Gut]. In fact,
 only the case where
  $s$ is large has been discussed there. In this case
 the most interesting phenomena  don't appear.
 For a while such problems were more or less forgotten.
During the last years some of these problems attracted interest
 again (see [NO], [O{\O}], [{\O}Z],[Dij]).

    I would like to  try to explain
 why this problem is interesting and also to
   discuss some approaches to this problem.

 {\bf 4.2. Limit as $s\rightarrow\infty.$}
  We want to show that the limit of the kernel-representations
   as $s\rightarrow +\infty$ is the
  canonical representation of $U(p,q)$
    in the space $L^2$ on Riemannian
   symmetric space $U(p,q)/(U(p)\times U(q))$.

     Consider the system of vectors
 $\Psi_z\in H_s$ (see Subsection 2.2). Let $\chi$ be a distribution
  in $B_{p,q}$ with compact support.    Consider the vector
    $\Theta(\chi)\in  H_s$
     defined by the equality
  $$
  \Theta(\chi)=\int_{B_{p,q}}
  |{\det}^{s}(1-zz^\ast)| \chi(z)\Psi_z  dz d\overline z
  $$
  and the scalar product $\{\cdot,\cdot\}_s$
  in the space of distributions on $B_{p,q}$
  with compact support given by the formula
\begin{eqnarray*}
& & \{\chi_1,\chi_2\}_s
:=  \langle\Theta(\chi_1),\Theta(\chi_2)\rangle  \\
&= & \int_{B_{p,q}}\int_{B_{p,q}}\Big| \frac
   {\det(1-zz^\ast)\det(1-uu^\ast)}
 {{\det}^2(1-zu^\ast)}
   \Big|^s \chi_1(z)\chi_2(u)
 dz d\overline z  du  d\overline u.
\end{eqnarray*}
   We can identify the space $H_s$ with
 the completion of the space   of distributions with respect
  to the scalar product $\{\cdot,\cdot\}_s$.
 The group $U(p,q)$ acts on this space
  of distributions by the formula
  $$
  B_s(g)f(z)=f(z^{[g]})
  $$
  (the formula doesn't depend on $s$,
 nevertheless, the scalar product
  and the spectrum of the representation
 do depend on $s$ in an essential way).

It appeared that for some normalizing multiplier
$\omega(s)$ for all continuous functions $\phi_1(z),\phi_2(z)$ on
$B_{p,q}$ with compact support we have
  $$
  \lim_{s\rightarrow +\infty}\omega(s)
  \{\phi_1,\phi_2\}_s
  =\int_{B_{p,q}}\phi_1(z)\phi_2(z)dzd\overline z
  $$

{\bf 4.3. Another formulation of the problem.}
Consider a spherical highest weight representation $T_s$
of the group $U(p,q)$ (see Subsection 2.3). Denote by $V_s$
the space of the representation $T_s$.  Denote by $T_s^\ast$ the
contragradient representation  of $T_s$.
Consider the tensor product $T_s\otimes T_s^\ast$.
This representation acts on the space $V_s\otimes V_s$ of
holomorphic functions on $B_{p,q}\times B_{p,q}$
via the operators
  $$
(T_s\otimes T_s^\ast) (g)
   = f(z_1,z_2)=f(z_1^{[g]},z_2^{[\overline g]})
   (a+z_1c)^{-s}(\overline a +z_2\overline c)^{-s}.
   $$
  The  $U(p,q)$-invariant scalar product in the space
  of holomorphic functions    on $B_{p,q}\times B_{p,q}$ is
   defined by the reproducing kernel
   $$
  M(z_1,z_2;u_1,u_2)=
  {\det}^{-s}(1-z_1u_1^\ast){\det}^{-s}(1-z_2 u_2^\ast)
  $$
  Consider the operator
  $$
  I:V_s\otimes  V_s\rightarrow H_s
  $$
  defined by the formula
  $$
  If(z)=f(z,\overline z)
  $$
  Obviously $I$ is a unitary operator interwining the representations
  $$T_s\ \otimes T_s^\ast \leftrightarrow A_s.$$

  Hence we can formulate our problem in the form:

 {\bf Problem $\bf 2'$.} {\it Decompose the tensor product $T_s\otimes T_s^\ast$}.

 {\bf 4.4. Orbits of the group
 $U(p,q)$ on the Shilov boundary.}

  Denote by $M_{p,q}$  the Shilov boundary of $B_{p,q}$.
  The elements of $M_{p,q}$ are matrices $z$ satisfying
the condition
    $$  z\cdot z^\ast=1.  $$
  In the other words, $z$ is  the matrix
 of an isometric embedding ${\mathbb C}^p\rightarrow {\C}^q$.
 Hence $M_{p,q}$ is a {\it complex Stiefel manifold}.

The Shilov boundary of $B_{p,q}\times B_{p,q}$
is $M_{p,q}\times M_{p,q}$. The group $U(p,q)$ has $(p+1)$ orbits
on $M_{p,q}\times M_{p,q}$.
 The only invariant of an orbit is the number
   \begin{equation}
   \alpha=rk(z-\overline u)
   \end{equation}
i.e if $(z_1,u_1),(z_2,u_2)\in M_{p,q}\times M_{p,q}$
and
$rk(z_1-\overline u_1)=rk(z_2-\overline u_2),$
then the points $(z_1,u_1),(z_2,u_2)$ belong to the same orbit.

  We denote by $\Xi_\alpha$ the orbit corresponding to a given
  invariant $\alpha$.

{\bf Remark.}    The orbit $\Xi_0$ is compact, whereas the
   orbit $\Xi_p$ is open (in the Shilov boundary).
For all $\alpha$  the closure of $\Xi_\alpha$
 is $\cup_{\sigma\le\alpha} \Xi_\sigma$.

{\bf 4.5. Restriction of holomorphic functions to $U(p,q)$-orbits
in the Shilov boundary.}
 Fix an orbit $\Xi_\alpha$ of   $U(p,q)$ in the Shilov boundary of
 $B_{p,q}\times B_{p,q}$. It can happen
   (and it really does happen) that for small
 $s$ the function
  $f\in H_s=V_s\otimes V_s$
    has  a well defined restriction to the orbit
   $\Xi_\alpha$. In this case the restriction operator
 is an interwining operator from $V_s\times V_s$
   to some Hilbert space of functions on $\Xi_\alpha$.

 It can also happen  (and it really does happen)
 that for small $s$ all  first    partial derivatives of
 the function $f\in H_s=V_s\otimes V_s$
   have a well defined restriction
  to the orbit $\Xi_\alpha$ etc. (see discussion
   of this phenomenon in [NO], Section 7).

   Fix $s$. For each $\alpha=0,1,\dots,p-1$
 consider the maximal number
   $\tau_\alpha$ such that all partial derivatives of order
$\le\tau_\alpha$   of functions $f\in V_s\otimes V_s$ have
 well-defined restrictions to $\Xi_\alpha$.
  (these numbers aren't known, but Theorem 3.1
  gives a possibility to estimate them;
   I don't know whether these estimates are best possible or not).
If restrictions of the functions $f\in H_s=V_s\otimes V_s$ to $\Xi_\alpha$
don't exist we suppose $\tau_\alpha=-1$.

{\bf Remark.}  For large $s$ the restriction operators don't
  exist, i.e. we have $\tau_\alpha=-1$ for all $\alpha<p-1$
(see below subsection 4.6).

Fix $\alpha=0,1,\dots\,p-1$. Consider some
$i_\alpha=1,\cdots,\tau_\alpha$.
 Denote by $Q[\alpha,i_\alpha]$
  the space of functions $f\in V_s\otimes V_s$
 such  that all partial   derivatives of $f$  of orders $\le i$ equal
zero on $\Xi_\alpha$.   We obtain a filtration
 \begin{equation}\label{filt}
\begin{array}{c}
  0\subset Q[p-1,\tau_{p-1}]\subset Q[p-1,\tau_{p-1}-1]\subset
    \dots\subset Q[p-1,0]\subset
      \\ \subset Q[p-2,\tau_{p-2}]\subset\dots\subset Q[p-2,1]\subset
 Q[p-2,0 ]\subset\dots
   \\
  \dots\subset Q[0,\tau_0]\subset\dots Q[0,1]\subset Q[0,0]
  \subset V_s\otimes V_s
    \end{array}\end{equation}

{\bf Remark.}   For large $s$ this filtration is trivial, but
  for small $s$ it is quite long.

  Consider the representations of $U(p,q)$
 in the subquotients of this filtration.
   Obviously $A_s=T_s\otimes T_s^\ast$
 is equivalent to the direct sum
   of the subquotients of the filtration (\ref{filt}).

{\bf Remark.}  The representations
  of $U(p,q)$ in  the subquotients
    have simple interpretations. For instance
   $(V_s\otimes V_s)/ Q[0,0]$
  is a subspace in the space
   of functions on the orbit $\Xi_0$.
 The space $Q[0,1]/Q[0,0]$ is a subspace in the space of sections
  of the normal bundle of the orbit $\Xi_0$ .
  The space $Q[0,1]/Q[0,2]$
    is a subspace in the  space of sections
     of the symmetric square of the normal bundle
       etc., see the discussion in [NO],\S 7.

  It is natural to hope that the spectrum
  in each subquotient is more or less
  `uniform', i.e. that the orbit structure of the
Shilov boundary    gives a separation of the
quite complicated spectrum of $A_s$ into
different   types (compare with the [GG]-project ).

  {\bf 4.6. Large $s$.} If $s$ is large
    enough, then the restriction
      operators don't exist. In this case the representation
  $A_s$ is equivalent to the standard
   representation of the group $U(p,q)$
   in $L^2$ on the Riemannian symmetric space
    $U(p,q)/(U(p)\times U(q))$, see [Ber2], [Rep], [Gut], [O{\O}].
   A sufficient (but not nessessary) condition for
  this is $s>p+q-1$ (i.e. $T_s$
   belongs to the Harish-Chandra discrete series).

    {\bf 4.6. Restriction to the compact orbit}. The part of the spectrum
  which corresponds to the compact orbit $\Xi_0$ is purely discrete
   and it consists of quite exotic representations of $U(p,q)$.
This part of the spectrum is not empty for small $s$ (for $s<(q-2p-1)/4)$.
In [NO] one can find a description of the analogous situation for
the compact orbit of $O(p,q)$ in the Shilov boundary of $B_{p,q}$.

In particular this gives   a  relatively simple way  to construct singular
unitary representations of the groups $U(p,q)$, see [NO].

 {\bf 4.7. Spherical kernel-representations of other classical groups}
  All classical Riemannian noncompact symmetric
   spaces $G/K$ (up to the center of $G$)  can be realized as
   (real, complex or quaternionic)  matrix balls (see [Ner3]).  Namely
  the space $G/K$ is the space
  of matrices $z$  over the field
   $\K=\R,  {\C},\H$ (see below) satisfying
    the additional condition (see below)      such that $\|z\|<1$.

     We list the symmetric spaces  $G/K$, fields and additional  conditions.

    $1^\ast.$  $U(p,q)/(U(p)\times U(q))$
       is the space of $p\times q$-matrices over ${\C}$.

   $2^\ast$ . $Sp(2n,{\R})/U(n)$
 is the space of symmetric
   $n\times n$-matrices over ${\C}$.

    $3^\ast$. $SO^\ast(2n)/U(n)$ is the space of
    skew symmetric $n\times n$-matrices over ${\C}$.

     $4^\ast$. $ O(p,q)/(O(p)\times O(q))$ is the space of
       $p\times q$-matrices over $\R$.

    $5^\ast$. $GL(n,{\R})/O(n)$
    is the space of symmetric $n\times n$-matrices
   over $\R$.

     $6^\ast$. $O(n,{\C})/O(n)$ is the space of skew-symmetric
       $n\times n$-matrices over $\R$.

   $7^\ast$. $GL(n,{\C})/U(n)$ is the space of hermitian
     $n\times n$-matrices over $\C$.

       $8^\ast$. $Sp(p,q)/(Sp(p)\times Sp(q))$ is the space of
    $p\times q$-matrices over $\H$.

     $9^\ast$. $GL(n,{\H})/Sp(p,q)$ is the space of  hermitian
       $n\times n$-matrices over  $\H$.

     $10^\ast$. $Sp(2n,{\C})/Sp(n)$ is the space of
      skew hermitian (i.e. $z=-z^\ast$)
       $n\times n$-matrices over $\H$.

    In all cases  enumerated above
     the group $G$ acts on the matrix ball
       $G/K$ by  fractional-linear transformations
         \begin{equation}  \label{frac}
           z\mapsto z^{[g]}:= (a+zc)^{-1}(b+zd)
             \end{equation}
  Now let us consider positive
  definite kernels on $G/K$ having the form
         $$
        L_s(z,u)=|\det (1-zu^\ast)|^{-2s}
               $$
  (the conditions for positive-definiteness are
   different for different spaces).
         Consider the Hilbert space defined
   by the positive definite kernel
    $L_s(z,u)$. We identify this space with a space of
    real analytic functions
     on $G/K$ in the same way as in 2.1.
      The group $G$ acts on $H_s$ by the unitary operators
             $$
 A_s(g)f(z)=f(z^{[g]})|\det(a+zc)|^{-2s}
    $$
     We call $A_s$ a {\it spherical kernel-representation of $G$}.

    {\bf 4.8. Nonspherical kernel-representations}.

    Consider one of the symmetric spaces form
     $2^\ast,3^\ast,5^\ast,6^\ast
     ,7^\ast,9^\ast,$ or $10^\ast$ and fix a finite dimensional
  irreducible representation
   $\rho $ of the group $GL(n,\K)$
   or of its universal covering.
    Assume that the function
    $$
    L(z,u)=\rho(1-zu^\ast)
    $$
    is a matrix valued positive definite kernel.
  Then we consider the associated Hilbert space
   $H_\rho$ of real-analytic functions
   $G/K\rightarrow \cal Y$
   and the unitary kernel-representation of $G$
    in $H_\rho$ given by the formula
    \begin{equation}\label{vectorr}
T_\rho (g)f(z)=\rho(a+zc)f(z^{[g]})
\end{equation}

Consider the cases $1^\ast,4^\ast,8^\ast$ and fix
a finite-dimensional  irreducible representation
 $\rho=\rho_1\otimes\rho_2$ of the
  group $GL(p,{\K})\times GL(q,\K)$ or
  of its universal covering.  Assume that the function
 $$
 L_{\rho}(z,u):=\rho_1(1-zu^\ast)\otimes\rho_2(1-u^\ast z)
 $$
 is a positive definite matrix-valued   kernel on  $G/K$.
  Then the group $G$ acts on the
   associated space of real-analytic
functions on $G/K$ by the formula
$$
T_\rho(g)f(z)=
\big(\rho_1(a+zc)\otimes\rho_2(d-cz^{[g]})\big)f(z^{[g]}).
$$

{\bf Remark.}  Our arguments from the Subsections 4.3-4.4
     are valid for general kernel-representations.

{\bf 4.9. Another description of kernel-representations.} (see [O{\O}])
For $G=Sp(2n,\R)$, $U(p,q)$, $SO^\ast(2n)$
 a kernel-representation is a tensor  product of an irreducible highest
  weight representation of $G$ and an  irreducible lowest weight
   representation of $G$.

In other cases a kernel-representation of $G$ is a restriction
 of a highest weight representation
  of the group $G^\ast$ to the {\bf symmetric} subgroup $G$:
$$\begin{array}{cccccc}
4^\ast. & G=O(p,q)  & G^\ast=U(p,q)
& 5^\ast.& G=GL(n,\R) &  G^\ast=Sp(2n,\R)\\
6^\ast.& G=O(n,\C) & G^\ast=SO^\ast(2n)&
7^\ast. &G=GL(n,\C)& G^\ast=U(n,n)\\
8^\ast. &G=Sp(p,q)&G^\ast=U(2p,2q)&
9^\ast. &G=GL(n,\H)&G^\ast=SO^\ast(2n)\\
10^\ast. &G=Sp(2n,\C)&G^\ast=Sp(4n,\R)&\ &\ &\
\end{array}$$

{\bf Remark.}   The cases $1^\ast-3^\ast$ can be described  in the same
 way.
 We have $G^\ast=G\times G$ and the embedding $G\rightarrow G^\ast$
   given by the formula $g\mapsto (g,g^\theta)$, where $\theta$ is
  an outer automorphism of $G$.

{\bf Remark.}  There are some
  additional possibilities related to
  highest weight representations of
   $O(p,2)$ and two exceptional groups
   (see [O{\O}]).

{\bf 4.10. The action of the Olshanskii semigroup}.
     For each matrix  ball
    $G/K$ consider the set of invertible matrices
 $g=\left(\begin{array}{cc}a&b\\c&d\end{array}\right)$
  such that the fractional-linear transformation
  (\ref{frac}) maps the matrix ball into itself.
   Obviously $\Gamma$ is a semigroup
 and the group $G$ is the group of invertible elements of $\Gamma$.
 The formula (\ref{vectorr}) defines a representation
 of the semigroup $\Gamma$. This
 representation is irreducible and all
 irreducible representations
 of $\Gamma$ can be obtained in  this way (see [Ols1]). See [Ols1],
 [Ner3], [Ner4], Appendix A, for the explicit
 description of the semigroups  $\Gamma$.
 Moreover the kernel-representations extend to representations of
 some categories (see [Ner3], [Ner4], Appendix A).
 I don't know any applications of these   phenomena to harmonic analysis
  of kernel-representations of groups.

  {\bf 4.11. Bibliographical comments.}

  a) Let $T$ be a highest weight
   representation of $G=Sp(2n,\R)$, $SO^\ast(2n)$, $U(p,q),\dots$
    and $S$ be a lowest weight representation of  $G$.
     Assume that $T,S$ be elements of the
      Harish-Chandra discrete series. Then
 $T\otimes S$ is equivalent to
 a representation of $G$ induced
  from an irreducible representation of the maximal compact
subgroup $K\subset G$  (see [Ber2], [Rep], [Gut]).

 b) Spherical kernel-representations for large $s$ are equivalent to
canonical representation of $G$ in $L^2(G/K)$ (see [O{\O}]).

c) Discrete spectra associated to
 the compact orbit in the Shilov boundary
 were  investigated in [NO]  for the case $G=O(p,q)$. Analogous
 results are valid for $G=U(p,q)$, $Sp(p,q)$ ([Ols2], [NO], 7.12.).
 One method of separation
  of the discrete spectrum is discussed in [Ner1],
  [NO],7.1-7.8. I think that the
   restriction operator to the compact
 orbit doesn't exist for   $G\ne U(p,q), Sp(p,q), O(p,q)$.

  d) The spectra for spherical    kernel-representations
of $U(2,2)$ were obtained in    [{\O}Z2].
For $1<s<3/2$ the spectrum consists of
 two different pieces. One of the pieces coincides
  with the spectrum of $L^2\big(U(2,2)/(U(2)\times U(2))\big)$.
   The other piece is an integral of nontrivial representations.

{\bf Question.} It is natural to think
  (it is not proved) that this piece of the spectrum
 is associated to the noncompact $U(2,2)$-orbit
  in the Shilov boundary
  of $B_{2,2}\times B_{2,2}$. Is it true?

  If this is so, then it is the only known
  case where a  spectrum associated to the
  noncompact orbit can be found.
   It is natural to think that such spectra
 exist in various spectral problems
 (not only for kernel-representations).

 e) The Spectrum of  spherical kernel-representations
  of $SU(1,1)$ was obtained in [VGG]
and in the case of $U(p,1)$ it was obtained in [Dij]
(the construction of the discrete increments was given in [Ner1]).

  f) For a Plancherel formula for  tensor product of highest
  weight and lowest weight (single-valued)  representations of
$SL(2,\R)$, see [Mol5].
Berezin [Ber2] announced a theorem which is almost equivalent to the
Plancherel formula for spherical kernel-representation in the case when $s$
is {\it large}  and the symmetric space is hermitian. Proof is
published in [UU].

  g) Nonspherical kernel-representations
  have discrete spectrum which
  is not associated to the compact orbit in the Shilov boundary.
   Some possibilities to find it
   are contained in the following two sections.
   A way to find the Harish-Chandra discrete series
    increments using trace theorems  is proposed in [Ner2].

 h) Let $\rho$ be the same as in 2.10.
  Let $\rho_s(\gamma)={\det}^{-s}(\gamma)\rho(\gamma)$.
   Then the limit of $T_{\rho_s}$ as $s\rightarrow\infty$
   is the representation of $G$
    induced from a finite dimensional
    representation of  $K$
    (i.e a representation on the space of  sections of
a vector bundle on $G/K$).
The theory of such representations is
 more or less equivalent to  Harish-Chandra's  theory for $L^2(G)$.

  e) {\it Deformation quantization,} see [Ber2]. As we have seen in 4.2
$H_s=V_s\otimes~V_s$. Hence we can identify the space
 $H_s$  with the space $V_s\rightarrow V_s$
of Hilbert-Schmidt operators.
The multiplication of Hilbert-Schmidt operators induces
an associative operation in $H_s$. This operation can be considered as
a ``deformation quantization", see [Ber3].

 \section{ Dual pairs}

 {\bf 5.1. Harmonic representation of $Sp(2n,\R)$}.
 We recall the definition of  the {\it harmonic representation}
 $W_{2n}$ (= Weil representation =
  Segal-Shale-Weil representation =
  Friedrichs-Segal-Berezin-Shale-Weil
   representation = oscillator representation )
    of the group $Sp(2n,\R)$ (see [KV], [Ner4]
     for a discussion of this representation).

 Consider the Hilbert space $F_n$ of holomorphic (entire) functions
 on ${\C}^n$  associated to   the reproducing kernel
 $\exp(\sum z_j \overline u_j)$ (this is the {\it boson Fock space}). Let
 $
 \left( \begin{array}{cc}\Phi&\Psi\\ \overline\Psi &
 \overline\Phi\end{array}\right)
 \in Sp(2n,\R)
 $ (see the realization of $Sp(2n,\R)$ described in Subsection 2.4)
 The harmonic representation  is a unitary
projective representation of the group $Sp(2n,\R)$
 in the space $F_n$ defined by the {\it Berezin formula} (see [Ner4]):
\begin{eqnarray}
& & W\left( \begin{array}{cc}\Phi&\Psi\\ \overline\Psi &
 \overline\Phi\end{array}\right) f(z)  \\
&=&  \int_{\C^N} \exp\{(z, \overline u)
 \left(\begin{array}{cc} \overline\Psi \Phi^{-1} &\Phi^{t-1}\\
 \Phi^{-1} & -\Phi^{-1}\Psi\end{array}\right)
 \left(\begin{array}{c}z\\ \overline u\end{array}\right)\Big\}
 f(z)\exp(-|z|^2) dz d\overline z.
\end{eqnarray}

{\bf 5.2.  Noncompact dual pairs.}
     Consider the following subgroups
  in the symplectic group {\it (noncompact
  Howe dual pairs)}:
  $$
  \begin{array}{c}
  Sp(2k(p+q),{\R})\supset Sp(2k,{\R})\times O(p,q)
  \\ Sp(2(k+l)(p+q),{\R})\supset U(k,l)\times U(p,q) \\
  Sp(4k(p+q),{\R})\supset SO^\ast(2k)\times Sp(p,q)
  \end{array}
  $$

{\bf Remark.}   Let us describe the first embedding. Consider the space
 ${\R}^{2n}$ equipped with a nondegenerate skew-symmetric bilinear
 form. Consider the space ${\R}^{p+q}$ equipped with a
 symmetric bilinear form with inertia indices
 $p,q$. The tensor product of symmetric and skew-symmetric
forms is a skew-symmetric form in the space
 ${\R}^{2n}\otimes{\R }^{p+q}$
 and  we obtain the required embedding of groups.
$\square$

  Let us restrict the harmonic representation $W_{2N}$
of the respective large symplectic group to these subgroups
   and then let us consider  the restriction  to the subgroups
   $$
   O(p,q),Sp(2k,{\R}), U(k,l), U(p,q), SO^\ast(2k), Sp(p,q)
   $$
    It was proved in [How], [Ada], [Li]
 that the spectra of these restrictions have discrete increments.
     This construction is one of the standard ways to obtain singular
unitary representations of the groups
     $U(p,q),O(p,q),Sp(p,q)$.

{\bf Proposition 5.1.} {\rm(see [NO])} {\it Each representation of
      $G= O(p,q)$, $Sp(2k,\R)$, $ U(k,l)$, $U(p,q)$, $SO^\ast(2k)$, $Sp(p,q)$
  which occurs in the spectrum of a dual pair discretely
   (resp. weakly) occurs in the spectrum
   of some kernel-representation discretely (resp.\ weakly). }

Proof. This proposition is more or less obvious.
    Consider for instance the case $Sp(2k,{\R})\times O(p,q)$.
    The restriction of $W_{2k(p+q)}$ to the subgroup
    $$
    Sp(2k,{\R})\subset Sp(2k,{\R})\times O(p,q)\subset Sp(2k(p+q),\R)
    $$
is equivalent to the representation
\begin{equation}\label{tenp}
W_{2k}^{\otimes p}\otimes (W_{2k}^\ast)^{\otimes q}
\end{equation}
The first tensor factor is a direct
 sum of highest weight representations
 and the second tensor factor is a direct sum of
 lowest weight representations.
  Hence (\ref{tenp}) is a direct sum of kernel-representations of
the group $Sp(2k,{\R})$.

Consider the subgroup $O(p,q)\subset Sp(2k,{\R})\times O(p,q)$.
 Consider the following subgroups in $Sp(2k(p+q),\R)$:
 $$
 \begin{array}{ccccccc}
 Sp(2k(p+q),\R)&\supset&O(p,q)&\times&Sp(2k,\R)&\supset&O(p,q)\\
 \| &\ &\cap &\ &\cup&\ &\cap
 \\Sp(2k(p+q),\R)&\supset&U(p,q)&\times &U(k)&\supset& U(p,q)
 \end{array}
 $$
 The restriction of $W_{2k(p+q)}$
  to $U(p,q)$ is a direct sum of
  highest weight representations and
  hence the restriction of $W_{2k(p+q)}$
 to $O(p,q)$ is a direct sum of kernel-representations.
$\square$

 {\bf 5.3. Restriction to orbits.}

 Consider the Cartan domain $C_N=Sp(2N,{\R})/U(N)$
and the reproducing kernel
  $$
  K(z,u)={\det}^{-1/2}(1-zu^\ast)
  $$
  on $C_N$. Denote by $H$ the associated Hilbert space.
   The group
   $Sp(2N,\R)$ acts on $H$ by the unitary operators
   $$
W^+_{2N}(g)=f(z^{[g]}){\det}^{-1/2}(\Phi+z\overline \Psi)
$$
The representation $W^+_{2N}$ is one of
 two irreducible components of the
 representation $W_{2N}$ (the other irreducible component
 can be realized in some space of holomorphic
 vector-valued functions). Again we have a question
  about the restrictions of holomorphic functions
  to $O(p,q)\times Sp(2k,\R)$-orbits in the Shilov boundary of
  $C_N$. The orbit structure
  of the Shilov boundary in this case is very complicated.
  In any case there exists an orbit
  $$
  {\cal F}=Sp(2k,{\R})/U(n) \times O(p,q)/Q,
  $$
  where $Q$ is stabilizer of a maximal
   isotropic subspace in the pseudoeuclidean
   space ${\R}^{p+q}$. I can show that the restriction operator
    restricting to the orbit
 $\cal F$ exists and this observation gives a way to see
 a part of the discrete spectrum for the dual pair.
 It is interesting to calculate this part of the spectrum.

 Another question which seems interesting to me: is it possible
 to obtain
by such way
 some handy realizations of some
 Harish-Chandra discrete series representations
in this way ?

 \section{The $L^2$-space on Stiefel manifolds}.

 {\bf 6.1. Stiefel manifolds.}
  The  following 10 types of homogeneous spaces $G/Q$
will be called {\it Stiefel manifolds}:
  $$
  \begin{array}{cccc}
  1^\circ .& O(p,q)/O(p-t,q-s) &
  2^\circ. & U(p,q)/U(p-t,q-s) \\
  3^\circ. & Sp(p,q)/Sp(p-t,q-s)&
  4^\circ. & Sp(2n,{\R})/Sp(2(n-t),{\R})\\
  5^\circ. & Sp(2n,{\C})/Sp(2(n-t),\C)&
  6^\circ.  & O(n,{\C})/O(n-t,\C)\\
  7^\circ.  & SO^\ast(2n)/SO^\ast(2(n-t))&\ & \ \\
  \end{array}
  $$

   $8^\circ-10^\circ$. The spaces of all linear embeddings
   $$
   \begin{array}{ccc}
   {\R}^{n-t}\rightarrow {\R}^n &
   {\C}^{n-t}\rightarrow {\C}^n
   &{\H}^{n-t}\rightarrow {\H}^n
   \end{array}$$

     In the last 3 cases the group $G$
is $GL(n,\R)$,$GL(n,\C)$,$GL(n,\H)$,
respectively, and $Q$ is the group of matrices having the form
\begin{equation}\label{q}
\left(\begin{array}{cc}
     1_s&\ast\\ 0 &\ast
\end{array}\right)
\end{equation}

{\bf Remark.}  The Stiefel manifold  $Sp(2n,{\R})/Sp(2(n-t),{\R})$
 is the space of isometric (form preserving) embeddings
 of the space ${\R}^{2t}$ equipped with
 a nondegenerate skew symmetric bilinear form on
 the space ${\R}^{2n}$
 equipped with a nondegenerated skew symmetric bilinear form. The
other  Stiefel manifolds $1^\circ-7^\circ$ have an analogous description.

 {\bf 6.2. Additional symmetries.}

We consider the homogeneous space $G/Q= Sp(2n,{\R})/Sp(2(n-t),{\R})$.
 Then the group $Sp(2t,\R)$ acts in the obvious
 way on the space  of symplectic-isometric embeddings ${\R}^{2t}
 \rightarrow {\R}^{2n}$ (since it acts on the space ${\R}^{2t}$).
 Hence the manifold\\  $ Sp(2n,{\R})/Sp(2(n-t),{\R})$
  is a $Sp(2t,{\R})\times Sp(2n,{\R})$-homogeneous space:
  $$
 Sp(2n,{\R})/Sp(2(n-t),{\R})=(Sp(2t,{\R})\times Sp(2n,{\R})){ /}
 ( Sp(2t,{\R})\times Sp(2(n-t),\R))
 $$
  Analogous additional groups of symmetries exist in all
the cases   $1^\circ-10^\circ$.
   These additional symmetries are useful since the spaces
   $L^2(G/Q)$ have $G$-spectrum with infinite multiplicities.

   {\bf 6.3.  $L^2$-spectrum on Stiefel manifolds}.
 Little is known about the spectral decomposition of $L^2(G/Q)$.
 Nevertheless,  this problem is interesting.
  See [Sch] and [Kob1] for  Flensted-Jensen type
  constructions of discrete spectra in $L^2$ on
  $$
  \begin{array}{ccc}
  O(p,q)/O(p-r,q), & U(p,q)/U(p-r,q), &
  Sp(p,q)/Sp(p-r,q)
  \end{array}
  $$
  For the case $r=1$ the Plancherel formula is obtained
  in [{\O}Z1].
  Some constructions  for discrete increments in spectra of
  $$
  { L^2(} U(p,q)/(U(p-t,q-t)\times U(p)\times U(q)){ )}
  \subset { L^2(} U(p,q)/(U(p)\times U(q)){ )}
  $$
  are contained in [RSW].

{\bf Remark.}  The cases  $G/Q$, where
$G=GL(n,\R)$, $GL(n,\C)$, $GL(n,\H)$
   are very simple. We can consider two-step unitary
induction from the subgroup $Q$ (see (\ref{q})) to
a parabolic subgroup
and from a parabolic to the group $GL(n)$ itself.

{\bf Proposition 6.1.} {\it Each representation of $G$ which is contained
 in the spectrum of $L^2(G/Q)$ discretely
 (resp. weakly) is contained in the spectrum
 of some kernel-representation of $G$ discretely (resp.\ weakly).}

Proof.  We use arguments from [How], [NO].
 Consider for instance the case $G=Sp(2n,\R)$.
  The representation
  $$
  W_{2n}\otimes W_{2n}^\ast
  $$
  of $Sp(2n,\R)$ is equivalent
 to the representation of $Sp(2n,\R)$ in
 $L^2({\R}^{2n})$. Hence the representation
 $$
 (W_{2n}\otimes W_{2n}^\ast)^{\otimes 2k}=
 W_{2n}^{\otimes 2k}\otimes (W_{2n}^\ast)^{\otimes 2k}
 $$
 is equivalent to the representation of
 $Sp(2n,\R)$ in $L^2$ on the space
 $Mat_{2k,2n}$ of all $2k\times 2n$-matrices.
 A generic orbit of $Sp(2n,\R)$
 in $Mat_{2k,2n}$ is a Stiefel
  manifold $Sp(2n,{\R})/Sp(2(n-k),\R)$.

{\bf Remark.}  For other groups $G$
   Proposition 6.1 can be proved using the same
 arguments. The basic observation is
 $$
 W_{2n}\Big|_{GL(n,{\R})}\simeq L^2({\R}^n)  $$
 (cf. the  real model of the harmonic representation in [KV]).

 {\bf 6.3. Some pseudo-Riemannian symmetric spaces.}
  In the obvious way we obtain
  $$
  \begin{array}{ccc}
  L^2\big( O(p,q)/(O(p)\times O(p-r,q))\big) &\subset &
   L^2\big( O(p,q)/O(p-r,q)\big)\\
   L^2\big(U(p,q)/(U(r)\times U(p-r,q))\big)&\subset
   & L^2\big( U(p,q)/U(p-r,q)\big)
   \\ L^2\big( Sp(p,q)/(Sp(r)\times Sp(p-r,q))\big)&\subset
   & L^2\big( Sp(p,q)/Sp(p-r,q)\big)
   \end{array}
   $$
   and the spectra of these spaces are
 contained in the spectra of Stiefel manifolds.
 I don't know such embeddings of
  spectra for other pseudo-Riemannian
  symmetric spaces.

\section{Miscellaneous}

{\bf 7.1. Molchanov representations.}
In 1970 Molchanov [Mol2] investigated the most degenerate representations
of the groups $O(p,q)$. These representations
are obtained by inducing from maximal  parabolic
subgroups and  taking unitary subquotients
of these induced representations. Molchanov representations
have very simple $K=O(p)\times O(q)$-spectra
(spherical harmonics with multiplicities less or equal $ 1$).
For  this reason Molchanov representations are very simple
and handy objects.

These representations are also very singular in different sences:
they have small rank in the sence of Howe (on Howe rank see [Li])
and their matrix coeffitients decrease very slowly (on decreasing of
matrix elements see [CHH]).

Recently these representations attracted
interest again, see [HT]. The most degenerate Molchanov
representations  are now very popular objects,
see  [BK], [Kos], [K{\O}]
(they are minimal representations in the sense of A.\ Joseph).
See in [BK] discussion of minimal representations of $O(p,q)$
and some other series of of simple groups.
The most peculiar case is  the `semi-exceptional' group $O(4,4)$,
see [Kos].

  I only wanted to mention that Molchanov representations
have very interesting spectral properties. Restrictions
of Molchanov representations to the subgroup $O(p,q-1)$
were discussed in [Mol] and [K{\O}]. The problem is
 solved for continuous series of Molchanov representations
(in this case we have continuous spectrum with
discrete increments, see [Mol2])
and for the most degenerate representation
(in this case the spectrum is discrete, see [K{\O}]).

  {\bf Questions.} It seems to me that it is a
very interesting  problem
to decompose the  tensor product of two
 Molchanov representations. Another interesting
question is `hunting' for discrete spectra
in tensor products of several Molchanov
representations
(paper [Ner1] contains a false remark on this subject).

{\bf Remark.}  It is possible to show that spectra of tensor
products of Molchanov representations contains spectra
of dual pairs. I don't know relations between
spectra of kernel-representations and spectra of tensor products
of Molchanov representations.

{\bf 7.2. Representations unitary with respect to indefinite scalar product.}
There exists a old folk-lor (originated by Ismagilov and Molchanov)
on harmonic analysis of represemtations which are unitary in
indefinite hermitian scalar products (partially it is published,
see [Ism], [Sul]). It is known that these representations often
(more often than on definite unitary case) contains discrete increments.
It was  observed in [Ner1] that this discrete increments usually
can be explained by restriction theorems. There arises very much
nonsolved problems (for instance it is possible to consider
kernel-representations with small or negative parameter $s$).
Neverless the main (well-known) question now
is the following:

{\bf Question.} Let a representation $\rho$  of a semisimple group
$G$ unitary in a indefinite hermitian metric is given.
When indefinite scalar product defines by some natural way a
topology in space of representation?

The question is very dangerous and sophisticated. On harmonic
analysis in handy
case of Pontryagin spaces there are several publications,
see for instance [Ism], [Sul] (a Pontryagin space is a linear space provided
by nondegenerate hermitian form with finite negative inertia index).
For some discussion of non-Pontraygin case (for highest weight representations)
see  [HN].

{\bf 7.3. Harish-Chandra  discrete series increments.} Tools which are proposed
in [NO] don't give possibilities to search Harish-Chandra discrete series
increments to unitary representations. Neverless some possibilities to expose
Harish-Chandra discrete series increments with restriction theorems
exist (this is the subject of [Ner2]).

{\bf 7.4. On totaly discrete spectra.} In [Ols4], [Kob2] there were discovered
large classes or restriction problems which have totaly discrete spectra.
Relations between these works are not clear to me, probably they
have open nondense intersection. It seems to me that it is interesting to
understand are these spectra related to some functional-theretical
phenomena or not.

\end{document}